\documentclass[journal,12pt,onecolumn,draftcls]{IEEEtran}
\usepackage{graphicx}
\usepackage{tikz}
\usepackage{amssymb,amsmath,cases}
\usepackage{algorithm}
\usepackage{algpseudocode}
\usepackage{algorithmicx}
\usepackage{amsfonts}
\usepackage[hidelinks]{hyperref}
\usepackage[normalem]{ulem}
\usepackage{cite}
\usepackage{enumitem}
\usepackage{subfig}
\usepackage{hyphenat}
\newcommand{\EE}{\mathcal{E}}
\newcommand{\GG}{\mathcal{G}}
\newcommand{\cA}{\mathcal{A}}
\newcommand{\cL}{\mathcal{L}}

\newcommand{\NN}{\mathcal{N}}

\newcommand{\cR}{\mathcal{R}}

\newcommand{\cS}{\mathcal{S}}

\newcommand{\R}{\mathbb{R}}
\newcommand{\diag}{\text{diag}}

\newcommand{\until}[1]{\{1,\ldots,#1\}}

\newcommand{\col}{\textsc{col}}
\newcommand{\blkdiag}{\textsc{blkdiag}}

\newcommand{\x}{x}
\newcommand{\w}{w}
\newcommand{\z}{z}
\newcommand{\s}{s}
\newcommand{\uu}{u}
\newcommand{\X}{X}

\newcommand{\uui}{\uu_{i}}
\newcommand{\dotuu}{\dot{\uu}}
\newcommand{\dotuui}{\dot{\uu}_i}

\newcommand{\phii}{\phi_{i}}

\newcommand{\1}{\mathbf{1}}

\newcommand{\norm}[1]{\left \|#1 \right \|}

\newcommand{\T}{^\top}
\newcommand{\RO}{\mathcal{R}_1}
\newcommand{\OR}{\mathcal{R}_2}

\newcommand{\alt}{z_{\text{alt}}}

\newcommand{\Go}{G_1}
\newcommand{\Gt}{G_2}

\newcommand{\Fsh}{F_{\sigma,h}}
\newcommand{\PL}{\mu}

\newcommand{\subj}{\textnormal{subj.~to}}

\newcommand{\tGt}{\tilde{G}_2}

\newcommand{\cX}{\mathcal{X}}

\newcommand\oprocendsymbol{\hbox{$\square$}}
\newcommand\oprocend{\relax\ifmmode\else\unskip\hfill\fi\oprocendsymbol}

\graphicspath{{figs/}}

\def\er/{Erd\H{o}s-R\'enyi}

\def\algo/{{\scshape Aggregative Tracking Feedback}}

\newtheorem{theorem}{Theorem}[section]

\newtheorem{lemma}[theorem]{Lemma}

\newtheorem{assumption}[theorem]{Assumption}
\newtheorem{remark}[theorem]{Remark}

\begin{document}

\title{Nonconvex Distributed Feedback Optimization\\ for Aggregative Cooperative Robotics}%
	
\author{Guido Carnevale, Nicola Mimmo, Giuseppe~Notarstefano %
\thanks{
	 This work was supported in part by the Italian Ministry of Foreign Affairs and International Cooperation, grant number BR22GR01.
Corresponding author: G.~Carnevale \texttt{guido.carnevale@unibo.it}.}
\thanks{G. Carnevale, N. Mimmo, and G. Notarstefano are with the Department of Electrical,  Electronic and Information Engineering,  Alma Mater Studiorum - Universita` di Bologna,  Bologna, Italy ({\tt \{guido.carnevale, nicola.mimmo2, giuseppe.notarstefano\}@unibo.it}).%
}}
\maketitle
  \begin{abstract}
    Distributed aggregative optimization is a recently emerged framework in
    which the agents of a network want to minimize the sum of local objective
    functions, each one depending on the agent decision variable (e.g., the local position of a team of robots) and an aggregation of all the agents'
    variables (e.g., the team barycentre).
    In this paper, we address a distributed feedback optimization framework in
    which agents implement a local (distributed) policy to reach a steady-state
    minimizing an aggregative cost function.
    We propose \algo/, i.e., a novel distributed feedback optimization law in
    which each agent combines a closed-loop gradient flow with a consensus-based
    dynamic compensator reconstructing the missing global information.
    By using tools from system theory, we prove that \algo/ steers the network
    to a stationary point of an aggregative optimization problem with (possibly)
    nonconvex objective function.
    The effectiveness of the proposed method is validated through numerical
    simulations on a multi-robot surveillance scenario.
  \end{abstract}

\section{Introduction}

The massive presence of networked systems in many areas is making distributed optimization more and more attractive for a wide range of tasks.
These tasks often involve dynamical systems (e.g., teams of robots or electric grids) that need to be controlled while optimizing a cost index.
Motivated by these scenarios, we propose a distributed feedback optimization law for aggregative problems.
Thus, we organize the literature in two parts involving (i) feedback optimization
and (ii) aggregative optimization.
Feedback optimization deals with an emerging class of controllers aiming at steering dynamic systems toward steady-states corresponding to solutions of associated optimization problems.
The key feature of these controllers is that they only rely on real-time gradient measurements, thus avoiding the knowledge of the cost functions and making them attractive in many fields such as control of electrical flow~\cite{dall2016optimal,tang2017real} and communication congestions~\cite{low2002internet}.
Early attempts designed these controllers through extremum-seeking techniques~\cite{wittenmark1995adaptive,teel2001solving,ariyur2003real,krstic2000stability,tan2006non}.
Feedback optimization laws for linear time-invariant systems are proposed in~\cite{menta2018stability,Lawrence2018Optimal,Colombino2020Online,Bianchin2022Time,cothren2022data}, 
while nonlinear plants are considered in~\cite{Jokic2009Constrained}.
The work~\cite{hauswirth2020timescale} studies the stability of different feedback optimization controllers based on gradient descent and its variations, including projected gradient and Newton method.
The works~\cite{haberle2020non,hauswirth2020anti,he2022model} also deal with the presence of constraints.
In~\cite{ospina2021data}, algebraic systems are controlled by relying on gradient information affected by random errors. %
As for feedback optimization in multi-agent systems, the early reference~\cite{Brunner2012Feedback} proposes an approach based on saddle point flows, while~\cite{terpin2022distributed} addresses a partition-based scenario. 
Despite cost-coupled (or consensus optimization) and
constraint-coupled setups are some of the most popular frameworks in distributed optimization~\cite{nedic2018distributed,yang2019survey,notarstefano2019distributed},
many applications in cooperative robotics cannot be suitably modeled in terms of these two setups.
Thus, the novel \emph{aggregative} optimization framework has recently gained attention. %
In this framework, agents in a network cooperatively
minimize the sum of local functions each depending on both a local variable and an \emph{aggregative} one.
The latter is obtained via aggregation of all the local variables (e.g., the mean). 
Differently from aggregative games, see, e.g.,~\cite{koshal2016distributed,gadjov2020single,
cenedese2020asynchronous,belgioioso2020distributed,carnevale2022tracking}, in the aggregative optimization the agents search for an optimal solution rather than a Nash equilibrium. 
This setup was introduced in~\cite{li2021distributed} in a static and unconstrained framework.
The static case is also studied in~\cite{chen2022distributed,wang2022distributed,grontas2022distributed}, while the online one is considered in~\cite{li2021distributedOnline,carnevale2021distributed,carnevale2022learning}.
The main contribution of the work is \algo/, i.e., a novel distributed
feedback optimization law for aggregative problems.
Its aim is to steer, in a fully distributed manner, a network of dynamic agents to a steady-state configuration which is a stationary point of a given aggregative optimization problem with nonconvex cost function.
In this scenario, each agent accesses local information only.
By relying on this information, our method (i) controls the network along an estimated descent direction of the cost, and (ii) reconstructs in each agent the global information needed for step (i).
Task (i) is performed through a distributed implementation of a closed-loop gradient flow. 
As per task (ii), a consensus-based dynamics asymptotically compensates for the mismatches between global and local information.
Differently from the existing
frameworks~\cite{li2021distributed,li2021distributedOnline,carnevale2021distributed,carnevale2022learning,chen2022distributed,wang2022distributed,grontas2022distributed}
in which agents do not have dynamics, \algo/ is a distributed feedback strategy handling at the same timescale the control and optimization of a network of nonlinear systems.
With tools from system theory, we guarantee the asymptotic
convergence of the network systems to a steady-state configuration corresponding to a stationary point of the problem.
Further, for isolated stationary points that are local minima, we prove asymptotic stability.
To the best of the authors' knowledge, this is the first work in the
literature proposing a distributed feedback law for a fully coupled
optimization problem.
Indeed, although~\cite{terpin2022distributed} provides a distributed optimization feedback law, a partition-based scenario is considered where the cost of each agent only depends on its own state and the neighboring ones.
Similar limitations appear in~\cite{Brunner2012Feedback}, which deals with convex problems and needs a central unit or a simplified problem structure.
Further, this is the first work dealing with aggregative optimization
problems with nonconvex cost functions.
Preliminary results related to this paper appeared
in~\cite{carnevale2022aggregative}, where, however, the proofs were omitted and a much simpler strongly convex setup with single-integrator dynamics was considered. 

The paper unfolds as follows. 
Section~\ref{sec:problem_formaulation} introduces the setup.
Section~\ref{sec:algo} presents
\algo/ and the main result of the work, while Section~\ref{sec:analysis} proves the latter with tools from system theory.
Finally, Section~\ref{sec:numerical_simulations} provides numerical simulations.\\ %

\paragraph*{Notation} $\col(v_1, \ldots, v_n)$ is the concatenation of the vectors $v_1, \ldots, v_n$.
$\diag(v)$ is the diagonal matrix with the components of $v$ on the diagonal.
$\blkdiag(M_1, \dots, M_N)$ is the matrix with $M_i \in \R^{n_i \times m_i}$ along the diagonal blocks. 
$\otimes$ is the Kronecker product. %
The identity matrix in $\R^{m\times m}$ is $I_m$.
$1_N$ and $0_N$ are the vectors of $N$ ones and zeros, respectively, while $\1 := 1_N \otimes I_d$.
Let $f: \R^{n_1}  \times \R^{n_2} \to \R^n$, then $\nabla_1 f(\x,y) := \dfrac{\partial}{\partial s}f(s,y)|_{s = x}$ and $\nabla_2 f(\x,y) := \dfrac{\partial}{\partial s}f(\x,s)|_{s = y}$. 
Let $X \subseteq \R^{n}$ and $x \in \mathbb{R}^n$, then $\|x\|_X := \inf_{y \in X}\norm{x - y}$.

\section{Problem Formulation and Preliminaries}
\label{sec:problem_formaulation}

This section describes the setup considered in this work.
We consider $N \in \mathbb{N}$ agents whose dynamics is
\begin{align}\label{eq:local_plant}
\dot{\x}_i = p_i(\x_i,\uui),
\end{align}
where $\x_i \in \R^{n_i}$ and $\uui \in
\R^{m_i}$ denote the state and the control input of agent $i$, while the dynamics $p_i:\R^{n_i} \times \R^{m_i}\to \R^{n_i}$ satisfies the following assumption.
\begin{assumption}\label{ass:steady_state}
        For all $i \in \until{N}$, there exist
        $h_i: \R^{m_i} \to \R^{n_i}$ and $\delta_1, \delta_2 > 0$ such that, for all $\bar{\uu}_i \in \R^{m_i}$, it holds $0 = p_i(h_i(\bar{\uu}_i),\bar{\uu}_i)$ and the trajectories of~\eqref{eq:local_plant} satisfy
		\begin{align*}
			\norm{\x_i(t) - h_i(\bar{\uu}_i)} \leq \delta_1\norm{\x_i(0) - h_i(\bar{\uu}_i)}\exp(-\delta_2t),
		\end{align*}
		for all $\x_i(0) \in \R^{n_i}$ and $t \ge 0$.		
		Further, there exist $L_h, L_p > 0$ such that
		\begin{align*}
		\|h_i(\uui)-h_i(\uui^\prime)\| &\le L_h \|\uui-\uui^\prime\|
		\\
		\|p_i(\x_i,\uui)-p_i(\x_i^\prime,\uui)\| &\le L_p \|\x_i - \x_i^\prime\|
		\\
		\|p_i(\x_i,\uui)-p_i(\x_i,\uui^\prime)\| &\le L_p \|\uui-\uui^\prime\|,
		\end{align*}
		for all $\x_i,\x_i^\prime \in \mathbb{R}^{n_i}$, $\uui,\uui^\prime \in \mathbb{R}^{m_i}$, and $i \in \until{N}$. 
\oprocend
\end{assumption}
Assumption~\ref{ass:steady_state} includes controlled systems embedding a regulator ensuring global exponential stability for the equilibrium $h_i(\bar{\uu}_i)$.
Thus, it is not a critical restriction.
The agents cooperate to reach a configuration representing a solution to the aggregative optimization problem
\begin{align}
	\begin{split}\label{eq:aggregative_problem}
	\min_{\substack{\x \in \R^n
	\\
	\uu \in \R^{m}}}  \: & \: \sum_{i=1}^{N}f_i(\x_i,\sigma(\x))
	\\
	\subj \: & \: \x_i = h_i(\uu_i), \forall i \in \until{N},
\end{split}
\end{align}
where, given $n := \sum_{i=1}^N n_i$, and $m := \sum_{i=1}^N m_i$, $\x := \col(\x_1, \dots, \x_N) \in \R^n$ and $\uu := \col(\uu_1,\dots,\uu_N)$ $\in \R^{m}$ are the global decision vectors, while $\sigma: \R^n \to \R^d$ is the \emph{aggregation function} defined as
\begin{align}\label{eq:sigma}
	\sigma(\x) = \dfrac{1}{N}\sum_{i=1}^{N}\phii(\x_i),
\end{align}
where $\phii : \R^{n_i} \to \R^d$ is the $i$-th contribution.
In order to lighten the notation, we introduce $F: \R^{n} \to \R$, $h: \R^{m} \to \R^{n}$, $\Fsh: \R^{m} \to \R$, and $G: \R^n \to \R^n$ defined as 
\begin{subequations}
	\begin{align}
		F(\x,\sigma(\x)) &:= \sum\nolimits_{i=1}^{N}f_i(\x_i,\sigma(\x))\label{eq:F}
		\\
		h(\uu) &:= \begin{bmatrix}
			h_1(\uu_1)\T& \dots& h_N(\uu_N)\T
		\end{bmatrix}\T
		\\
		\Fsh(\uu) &:= F(h(\uu),\sigma(h(\uu)))\label{eq:Fsh}
		\\
		G(\x) &:= \nabla F(v,\sigma(v))\mid_{v = \x}.
	\end{align}
\end{subequations}
According to the distributed computation paradigm, we assume that the global information of problem~\eqref{eq:aggregative_problem} is not locally available for the single agent $i$. 
We remark that, differently from the existing works~\cite{li2021distributed,li2021distributedOnline,carnevale2021distributed,carnevale2022learning,chen2022distributed,wang2022distributed,grontas2022distributed}, we also satisfy the feedback optimization paradigm in the following sense.
The analytic expression of the local objective functions and aggregation rules are not available to the agents, they can be only measured according to current local variables.
In particular, each agent $i$ can only access $\nabla_1 f_i (\x_i,\hat{\sigma}_i)$,
$\nabla_2 f_i(\x_i,\hat{\sigma}_i)$, $\phii(\x_i)$, $\nabla\phii (\x_i)$, and $\nabla h_i(\uu_i)$ where $\x_i$ and $\uu_i$ are the current state and input variables, respectively, while $\hat{\sigma}_i \in \R^d$ is its local estimate of $\sigma(\x)$.
\begin{assumption}\label{ass:lipschitz}
$\Fsh$ is radially unbounded and differentiable and there exist $L_0,L_1,L_2 > 0$ such that
\begin{align*}
		\|\nabla F(\x,\sigma(\x))  -  \nabla F(\x^\prime,\sigma(x^\prime))\|  &\le   L_0 \|x-x^\prime\|
		\\
\norm{\nabla_1 f_i(\x_i,y_i) - \nabla_1 f_i(\x_i^\prime,y_i^\prime)} &\le  L_1\norm{\begin{bmatrix}\x_i - \x_i^\prime\\
	 y_i - y_i^\prime\end{bmatrix}}
\\
\|\nabla_2 f_i(\x_i,y_i)-\nabla_2 f_i(x_i^\prime,y_i^\prime)\| &\le  L_2 \norm{\begin{bmatrix}\x_i - \x_i^\prime\\
	 y_i - y_i^\prime\end{bmatrix}},
\end{align*}
for all $\x, \x^\prime \in \R^n$, $y, y^\prime \in \R^{Nd}$, $\x_i, \x^\prime_i
\in \R^{n_i}$, $y_i, y^\prime_i \in \R^{d}$, and $i \in \until{N}$. 
Also, the functions $\phii$ are differentiable and there exists $L_3 > 0$ such that
	\begin{equation*}
	\norm{\phii(\x_i) - \phii(\x_i^\prime)} \leq L_3\norm{\x_i - \x_i^\prime},
	\end{equation*}
	for all $\x_i$, $\x_i^\prime \in \R^{n_i}$ and $i \in \until{N}$. \oprocend
\end{assumption}
Assumption~\ref{ass:lipschitz} imposes very mild requirements about the class of functions considered in this work.
Indeed, differently from~\cite{li2021distributed,li2021distributedOnline,carnevale2021distributed,carnevale2022learning,chen2022distributed,wang2022distributed,grontas2022distributed,carnevale2022aggregative}, the cost convexity is not required.
The enforced Lipschitz properties are the same.
The communication among the agents is performed according to a directed graph
$\GG = (\until{N}, \EE)$ with $\EE \subset \until{N}\times \until{N}$ being the
edge set. 
If an edge $(j,i)$ belongs to $\EE$, then agent $i$ can receive
information from agent $j$, otherwise not. 
The set of (in-)neighbors of agent $i$ is
defined as $\NN_i := \{j \in \until{N} \mid (j,i) \in \EE\}$. 
We associate to the graph $\GG$ a weighted adjacency matrix $\cA \in\R^{N\times N}$ whose
entries satisfy $a_{ij} >0$ whenever $(j,i)\in \EE$ and
$a_{ij} =0$ otherwise. 
The weighted in-degree and out-degree of agent $i$ are
defined as $d_i^{\text{in}} = \sum_{j\in \NN_i}a_{ij}$ and
$d_i^{\text{out}}= \sum_{j\in \NN_i}a_{ji}$, respectively. 
Finally, we associate to $\GG$ the so-called
Laplacian matrix defined as $\cL:=\mathcal{D}^{\text{in}} - \cA$, where
$\mathcal{D}^{\text{in}} := \diag(d_1^{\text{in}},\ldots,d_N^{\text{in}}) \in \R^{N \times N}$.
The following assumption about network connectivity details the class of graphs considered in this work. %
\begin{assumption}\label{ass:network}
	$\GG$ is strongly connected and weight-balanced, namely $d_i^{\text{in}} = d_i^{\text{out}}$ for all $i \in \until{N}$.
	\oprocend
\end{assumption}

\section{\algo/}
\label{sec:algo}

\begin{figure}[htpb]
	\centering
	\includegraphics[scale=1]{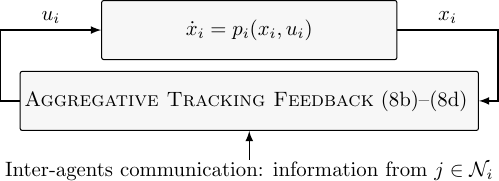}
	\caption{Block diagram describing~\eqref{eq:continuous_aggregative}.}
	\label{fig:local_scheme}
\end{figure}

This section describes \algo/, i.e., a distributed feedback optimization
law designed to steer the agents (with dynamics~\eqref{eq:local_plant}) to a configuration corresponding to a solution of problem~\eqref{eq:aggregative_problem}.
To introduce our law, given any $\uu_i \in \R^{m_i}$, let us study problem~\eqref{eq:aggregative_problem} when $\x_i = h_i(\uui)$ for all $i \in \until{N}$, i.e., when each agent has already reached its steady-state configuration (see Assumption~\ref{ass:steady_state}). 
Then,~\eqref{eq:aggregative_problem} becomes
\begin{align}\label{eq:tilde_aggregative_problem}
\min_{\uu \in \mathbb{R}^m} \sum_{i=1}^N f_i(h_i(\uui),\sigma(h(\uu))).
\end{align}
It is well-known that~\eqref{eq:tilde_aggregative_problem} can be addressed by adopting the continuous-time gradient method (see, e.g.,~\cite{bloch1994hamiltonian}), which, for all $i \in \until{N}$, reads as
\begin{align}\label{eq:centralized_aggregative_with_u}
\dotuui &= - \dfrac{\partial}{\partial \uui} F(h(\uu),\sigma(h(\uu)))
\notag\\
&= - \nabla h_i(\uui) \Big(\nabla_1 f_i(h_i(\uui),\sigma(h(\uu))) 
+\dfrac{\nabla \phii(h_i(\uui))}{N} \sum_{j=1}^{N}\nabla_2 f_j(h_j(\uu_j),\sigma(h(\uu))) 
\Big).
\end{align}
However, 
as customary in the context of feedback optimization, agent $i$ can access $\nabla h_i(\uui)$ but the quantities $\nabla_1 f_i(h_i(\uui),\sigma(h(\uu)))$, $\nabla \phii(h_i(\uui))$, and $\sum_{j=1}^{N}\nabla_2 f_j(h_j(\uu_j),\sigma(h(\uu)))$ are not available.
Indeed, agent $i$ can only \emph{measure} $\nabla_1 f_i(\x_i,\hat{\sigma}_i)$, $\nabla \phii(\x_i)$, and $\sum_{j=1}^{N}\nabla_2 f_j(\x_j,\hat{\sigma}_j)$, where $\hat{\sigma}_i \in \R^{d}$ is the local estimate of agent $i$ about $\sigma(\x)$.
Thus, we modify~\eqref{eq:centralized_aggregative_with_u} as
	\begin{align}
		\dotuui  &= - \nabla h_i(\uui)\Big(\nabla_1 f_i(\x_i,\hat{\sigma}_i)
		+\dfrac{\nabla \phii(\x_i)}{N} \sum_{j=1}^{N}\nabla_2 f_j(\x_j,\hat{\sigma}_j) 
		\Big).\label{eq:centralized_aggregative}
	\end{align}
The control law~\eqref{eq:centralized_aggregative} cannot be implemented in a distributed fashion because $\sum_{j=1}^{N}\nabla_{2}f_j(\x_j,\hat{\sigma}_j)$ represents a global information. %
Therefore, law~\eqref{eq:centralized_aggregative} needs to be properly intertwined with an estimation mechanism providing both $\sigma(\x)$ and $\sum_{j=1}^{N}\nabla_{2}f_j(\x_j,\hat{\sigma}_j)$.
To this end, we take inspiration from the continuous-time compensation dynamics proposed in~\cite{carnevale2023triggered} which locally reconstructs the unavailable global gradient in a consensus optimization framework.
Then, we embed two consensus-based mechanisms (with local states $\w_i, \z_i \in \R^{d}$) giving rise to the distributed feedback optimization law termed \algo/ and resumed in Algorithm~\ref{algo:algo} from the perspective of agent $i$.
The parameter $\alpha_1 >0$ tunes the speed of $\uu_i$ and $(\w_i,\z_i)$ relative to the plant.
In turn, the parameter $\alpha_2 > 0$ tunes the speed of $(\w_i,\z_i)$. %
As it will become clearer with the formal analysis, the role of $\alpha_1$ and $\alpha_2$ is to impose two timescale separations between the dynamics~\eqref{eq:local_plant},~\eqref{eq:x_i_dynamics}, and~\eqref{eq:w_i_dynamics}-\eqref{eq:z_i_dynamics}.
More in detail, the parameters $\alpha_1$ and $\alpha_2$ are tuned to make the dynamics~\eqref{eq:w_i_dynamics}-\eqref{eq:z_i_dynamics} the fastest ones, and the dynamics~\eqref{eq:x_i_dynamics} the slowest one.
The role of the initialization $w_i(0) = z_i (0) = 0_d$ for all $i \in \until{N}$ will be detailed into Section~\ref{sec:reformulation}.
The weights $a_{ij}$, with $j \in \mathcal{N}_i$, are the entries of the weighted adjacency matrix $\cA$ associated to the graph $\GG$ and, thus, belong to the available local information of agent $i$.
Agent $i$ exchanges with its neighbors the information $(w_i + \phii(\x_i))$ and $\z_i + \nabla_2 f_i(\x_i,\w_i + \phii(\x_i))$, i.e., a total of $2d$ components.
Although Algorithm~\ref{algo:algo} requires continuous-time inter-agents communication, we note that its algorithmic structure lends itself to extensions implementing discrete-time event-triggered communication, see, e.g.,~\cite{kia2015distributed,carnevale2023triggered}.
Fig.~\ref{fig:local_scheme} describes the closed-loop system~\eqref{eq:continuous_aggregative} in terms of block-diagrams.
\vspace{-.1cm}
\begin{algorithm}[htpb]
	\begin{algorithmic}
		\State initialization: 
		\State $\x_i(0) \in \R^{n_i}, \uu_i(0) \in \R^{m_i}$, $\w_i(0) = 0_d, \z_i(0) = 0_d$
	\begin{subequations}\label{eq:continuous_aggregative}
		\begin{align}
			\dot{\x}_i &= p_i(\x_i,\uui)\label{eq:state_dynamics}
			\\
				\dotuu_i  &= - \alpha_1\nabla h_i(\uui) \left(\nabla_1 f_i(\x_i,\w_i+\phii(\x_i))
			+\nabla\phii(\x_i) \left(\z_i+\nabla_2f_i(\x_i,\w_i+\phii(\x_i))\right)
			\right)\label{eq:x_i_dynamics}\\
			\dot{\w}_i &= -\frac{1}{\alpha_2}\sum_{j\in\NN_i}a_{ij}\left(\w_i + \phii(\x_i) - \w_j - \phi_j(\x_j)\right)
			\label{eq:w_i_dynamics}
			\\
			\dot{\z}_i &=
			-\frac{1}{\alpha_2}\sum_{j\in\NN_i}  a_{ij}(z_i  +  \nabla_{2}f_i(\x_i,\w_i+ \phii(\x_i)))
			+\frac{1}{\alpha_2}\sum_{j\in\NN_i}  a_{ij} (z_j  +  \nabla_{2}f_j(\x_j,\w_j+ \phi_j(\x_j)))\label{eq:z_i_dynamics}
		\end{align}
	\end{subequations}
	\end{algorithmic}
	\caption{Agent $i$ dynamics in closed-loop}
	\label{algo:algo}
\end{algorithm}
\vspace{-.27cm}
\vspace{-.27cm}

The next theorem provides the convergence and stability properties of the closed-loop system resulting from~\eqref{eq:continuous_aggregative}.
To this end, let $\mathcal{S} := \{\col(\x,\uu,\w,\z) \in \R^{n+m+2Nd}\mid  \1\T \w = 0, \1\T \z = 0  \}$, $X \subset \R^{n} \times \R^{m}$ 
be the set of stationary points of problem~\eqref{eq:aggregative_problem}, $\cX := \{\col(\x,\uu,\w,\z) \in \R^{n+m+2Nd} \mid (x,u) \in X, \col(w,z) = 0_{2Nd}\}$, and 
\begin{align*}%
	&\w := \begin{bmatrix}
		\w_1\\
		\vdots
		\\
		\w_N
\end{bmatrix}, \hspace{.1cm} \z := \begin{bmatrix}
		\z_1\\
		\vdots\\
		\z_N 
	\end{bmatrix}, \hspace{.1cm} \pi^w(x) := \begin{bmatrix}
		-\phi_1(x_1) + \sigma(x)
		\\
		\vdots 
		\\
		-\phi_N(x_N) + \sigma(x)
	\end{bmatrix}
	\\
	&\pi^z(x) := \begin{bmatrix}
		-\nabla_2 f_1(\x_1,\sigma(x)) + \frac{1}{N}\sum_{j=1}^{N}\nabla_2 f_j(\x_j,\sigma(x)) 
		\\
		\vdots
		\\
		-\nabla_2 f_N(\x_N,\sigma(x)) + \frac{1}{N}\sum_{j=1}^{N}\nabla_2 f_j(\x_j,\sigma(x)) 
	\end{bmatrix}.
\end{align*}
\begin{theorem}\label{th:convergence}
	Consider~\eqref{eq:continuous_aggregative} and let Assumptions~\ref{ass:steady_state},~\ref{ass:lipschitz}, and~\ref{ass:network} hold.
	Then, there exist $\bar{\alpha}_1, \bar{\alpha}_2 > 0$ such that, for all $\col(\x_i(0),\uu_i(0),\w_i(0),\z_i(0)) \in \R^{n_i + m_i + 2d}$ such that $\z_i = \w_i = 0$ for all $i \in \until{N}$, $\alpha_1 \in (0,\bar{\alpha}_1)$, and $\alpha_2\in(0,\bar{\alpha}_2)$, the trajectories are bounded and it holds
	\begin{align}\label{eq:limit_statement}
		\lim_{t\to\infty} \norm{\begin{bmatrix}
			\x(t)
			\\
			\uu(t)
			\\
			\w(t) - \pi^{\w}(x(t)) 
			\\
			\z(t) - \pi^{\z}(x(t))
		\end{bmatrix}}_{\cX}&= 0.
	\end{align}
	Further, given any $\bar{u} \in \R^{m}$ being an isolated stationary point and a local minimum of $\Fsh$, the point $\col(h(\bar{u}),\bar{u}, \pi^{\w}(h(\bar{u})),\pi^{\z}(\bar{u}))$ is locally asymptotically stable for the closed-loop system resulting from~\eqref{eq:continuous_aggregative} restricted to $\cS$.
\end{theorem}
The proof is given in Section~\ref{sec:proof}.
Theorem~\ref{th:convergence} guarantees that, under a suitable initialization, \algo/ asymptotically steers the network into the set $\X$ of stationary points of
problem~\eqref{eq:aggregative_problem}.
More in detail, with the initial conditions restricted to $\mathcal{S}$, Theorem~\ref{th:convergence} guarantees the global asymptotic stability of the set $\{(\x,\uu,\w,\z) \in \R^{n+m+2Nd} \mid (\x,\uu) \in X, \w = \w = \pi^{\w}(\x), \z = \pi^{\z}(\x)\}$~\cite[Th.~1]{Yuandan1994Recent} for the networked dynamics arising from~\eqref{eq:continuous_aggregative}. %
Further, with the initial conditions restricted to $\mathcal{S}$, Theorem~\ref{th:convergence} ensures asymptotic stability for the configurations corresponding to isolated stationary points $\bar{\uu}$ of $\Fsh$ that are local minima, i.e., the point $\col(h(\bar{u}),\bar{u}, \pi^{\w}(h(\bar{u})),\pi^{\z}(\bar{u}))$.
Namely, for each local minimizer $\bar{\uu}$, there exists $\mathcal{S}_0 \subseteq \mathcal{S}$ such that the system trajectories starting in $\mathcal{S}_0$ asymptotically converge to $\col(h(\bar{u}),\bar{u}, \pi^{\w}(h(\bar{u})),\pi^{\z}(\bar{u}))$.
\begin{remark}
	Theorem~\ref{th:convergence} 
	is in agreement with the discussions in~\cite{hauswirth2020timescale} and~\cite{absil2006stable}.
	More in detail, although any stationary point has a related configuration representing an equilibrium of the aggregate form of~\eqref{eq:continuous_aggregative}, in general, the configurations related to local minima are stable equilibria but not necessarily asymptotically stable.
	Further, the configurations related to saddle points and local maxima are unstable.
	Finally, in the case of (local) strong convexity, the configuration corresponding to the problem solution is exponentially stable.
	This last claim holds since the strong convexity in a neighborhood $\mathcal{U}$ of a minimizer $\bar{\uu} \in \R^m$ of $\Fsh$ implies that $\norm{\nabla \Fsh(\uu)}^2 \ge 2\PL(\Fsh(\uu) - \Fsh(\bar{\uu}))$ for all $\uu\in \mathcal{U}$ and some $\PL > 0$.
	By plugging this result into the Lyapunov function derivative bound provided in~\eqref{eq:Vdot_negative}, the exponential stability follows by~\cite[Th.~4.10]{khalil2002nonlinear}.\oprocend
\end{remark}
The key intuition behind Theorem~\ref{th:convergence} consists in reformulating the whole closed-loop dynamics resulting from~\eqref{eq:continuous_aggregative} as a singularly perturbed system given by the interconnection of three subsystems. %
In detail, the bounds $\bar{\alpha}_1$ and $\bar{\alpha}_2$ characterize the required timescale separation and depend on the problem parameters (e.g., network connectivity, Lipschitz constants of the dynamics, etc.).

\section{Closed-Loop System Analysis}
\label{sec:analysis}

In this section, we prove Theorem~\ref{th:convergence} by using tools from system theory. 
Assumptions~\ref{ass:steady_state},~\ref{ass:lipschitz}, and~\ref{ass:network} hold for the entire section.
The main steps of the proof are:

\begin{enumerate}[label=(\roman*)]
\item We rewrite the aggregate version arising from~\eqref{eq:continuous_aggregative} as the interconnection of
  three dynamic subsystems describing the evolution of all the variables.
\item Within three separate lemmas we give suitable properties of the
  time-derivative of three different Lyapunov-like functions. Specifically, each
  one of these lemmas assesses the stability properties 
	of the three subsystems
  identified within step (i). %
\item To conclude, we define a candidate Lyapunov function for the whole system
  and, relying on the lemmas of step (ii) and LaSalle arguments, we study its
  time-derivative to prove Theorem~\ref{th:convergence}.
\end{enumerate}
Steps (i) and (ii) are carried out in Section~\ref{sec:reformulation} and~\ref{sec:preparatory}. 
Section~\ref{sec:proof} is devoted to the development of step (iii).

\subsection{System Reformulation}
\label{sec:reformulation}

In this section, we reformulate~\eqref{eq:continuous_aggregative} by leveraging the initialization of $w$ and $z$ and the consensus properties of their dynamics.
To this end, we define $L = \cL \otimes I_{d}$ and the operators %
$\Go: \R^{n} \times \R^{Nd} \to \R^{n}$ and $\Gt: \R^{n} \times \R^{Nd} \to \R^{n}$ given by
	\begin{align*}
	\Go(\x,\s) &=
	\begin{bmatrix}
	\nabla_1 f_1(\x_1, \s_1)
	\\
	\vdots
	\\
	\nabla_1 f_N(x_N, \s_N)
\end{bmatrix}, \quad
	\Gt(\x,\s) = \begin{bmatrix}
		\nabla_2 f_1(\x_1,\s_1)
		\\
		\vdots\\
		\nabla_2 f_N(x_N, \s_N)
	\end{bmatrix},
\end{align*}
where $\x = \col(\x_1,\dots,\x_N)$ and $\s = \col(\s_1,\dots,\s_N)$ with $x_i \in \R^{n_i}$ and $\s_i \in \R^d$ for all $i \in \until{N}$.
Then, the stacked column form of~\eqref{eq:continuous_aggregative} reads as
\begin{subequations}\label{eq:continuous_aggregative_gloabal}
	\begin{align}
		\dot{\x} &= p(\x,\uu)\label{eq:x_dynamics_before_change}
			\\
			\dotuu &= -\alpha_1\nabla h(\uu)\left(\Go(\x,\w + \phi(\x)) + \nabla \phi(\x)\left(\z+\Gt(\x,\w + \phi(\x))\right)\right)\label{eq:u_dynamics_before_change}
			\\
			\dot{\w} &= -\frac{1}{\alpha_2}L\left(\w+\phi(\x)\right) \label{eq:w_dynamics}
			\\
			\dot{\z} &= -\frac{1}{\alpha_2}L\left(\z+\Gt(\x,\w+ \phi(\x))\right).\label{eq:z_dynamics}
	\end{align}
\end{subequations}
	We note that $\mathcal{S}$ is forward-invariant for~\eqref{eq:continuous_aggregative_gloabal} because $\1^\top L = 0$ (cf. Assumption~\ref{ass:network}). 
	Hence, let us take advantage of this property via a change of variables.
	To this end, let $R \in \R^{Nd \times (N-1)d}$ be such that $R\T R = I$ and
	$R\T \1 = 0$, and $T \in \R^{2Nd \times 2Nd}$ be defined as $T:= \begin{bmatrix}R& \1/N\end{bmatrix}\T$.
	Then, let $\eta, \zeta \in \mathbb{R}^{(N-1)d}$, $\eta_{\text{avg}}, \zeta_{\text{avg}} \in \mathbb{R}^{d}$ be
	\begin{align}\label{eq:mean_change_of_variables}
		\begin{bmatrix}
			\eta\\
			\eta_{\text{avg}}
		\end{bmatrix}
		:= 
		T\w ,
		\quad	
		\begin{bmatrix}
			\zeta\\
			\zeta_{\text{avg}}
		\end{bmatrix} := T\z. 
	\end{align}
	Then, by using~\eqref{eq:w_dynamics}-\eqref{eq:z_dynamics}, it holds $\dot{\eta}_{\text{avg}}= 0$ and $\dot{\zeta}_{\text{avg}} = 0$.
	Thus, the initialization $w(0)=z(0)=0_{Nd}$ guarantees that ${\eta}_{\text{avg}}(t)  =  \1 {\zeta}_{\text{avg}}(t)  =  0_d$ for all $t  \ge  0$.
	Then, by defining $\psi := \col(\eta,\zeta)$ and using~\eqref{eq:mean_change_of_variables}, we restrict~\eqref{eq:w_dynamics}-\eqref{eq:z_dynamics} to
			\begin{align}\label{eq:psi_dynamics}
				\dot{\psi} &= \frac{1}{\alpha_2}\begin{bmatrix}
					-R\T LR& 0\\
					0& -R\T LR
				\end{bmatrix}\psi
				+ \frac{1}{\alpha_2}\begin{bmatrix}
					-R\T L& 0\\
					0& -R\T L
				\end{bmatrix}\begin{bmatrix}
					\phi(\x)\\
					\Gt(\x,\begin{bmatrix} R& 0\end{bmatrix}\psi + \phi(\x))
				\end{bmatrix}.
			\end{align}
		Let $\cR := \blkdiag(R\T,R\T)$. 
		Then, we note that
		\begin{align}
			\label{eq:barpsi}
			\bar{\psi}(\x) := -\cR\col(\phi(\x),\Gt(\x,\1\sigma(\x)))
		\end{align}
		represents an equilibrium for \eqref{eq:psi_dynamics} for all $\x \in \R^n$. 
		Hence, let the error coordinate $\xi \in \R^{2(N-1)d}$ and $\RO, \OR \in \R^{Nd \times 2(N-1)d}$ be defined as
			$\xi := \psi - \bar{\psi}(\x)$
and $\RO := \begin{bmatrix}
	R& 0
\end{bmatrix}$ and $\OR := \begin{bmatrix}
	0& R
\end{bmatrix}$.
Then, by using this notation and
exploiting~\eqref{eq:sigma},~\eqref{eq:psi_dynamics}, and $I-R R^\top = \1 \1^\top / N$, we
equivalently rewrite~\eqref{eq:continuous_aggregative_gloabal} as
		\begin{subequations}\label{eq:continuous_aggregative_gloabal_man}
			\begin{align}
				\dot{\x} &= p(\x,\uu)\label{eq:xx_dynamics}
					\\
					\dotuu &= -\alpha_1\nabla h(\uu)\left(\Go(\x,\RO\xi + \1\sigma(\x)) 
					+ \nabla \phi(\x)\left(\dfrac{\1\1\T}{N}\Gt(\x,\1\sigma(\x))+\OR\xi + \tGt(\x,\xi) \right)\right)\label{eq:uu_dynamics}
					\\
				\dot{\xi} &= \frac{1}{\alpha_2}\begin{bmatrix}
					-R\T LR& 0\\
					0& -R\T LR
				\end{bmatrix}\xi
				+ \frac{1}{\alpha_2}
				\begin{bmatrix}
					0\\
					R\T L \tGt(\x,\xi)
				\end{bmatrix}
				-\nabla\bar{\psi}(\x) \,p(\x,\uu),\label{eq:xi_dynamics}
			\end{align}
		\end{subequations}
		where $\tGt: \R^{n} \times \R^{2(N-1)d} \to \R^{Nd}$ is defined as
		\begin{align*}
			\tGt(\x,\xi) := \Gt(\x,\RO\xi + \1\sigma(\x)) - \Gt(\x,\1\sigma(\x)).
		\end{align*}

\subsection{Preparatory Results}
\label{sec:preparatory}

In this section, we provide three preparatory results needed to prove Theorem~\ref{th:convergence}.
\begin{lemma}\label{lemma:W}
  There exists $W: \R^n \times \R^m \to \R$ such that, along the trajectories of~\eqref{eq:xx_dynamics} and~\eqref{eq:uu_dynamics}, it holds
  \begin{subequations}\label{eq:W_conditions_lemma}
    \begin{align}
      c_1 \norm{\x - h(\uu)}^2 &\leq W(\x,\uu) \leq c_2\norm{\x - h(\uu)}^2\label{eq:W_upepr_bound}
      \\
      \dot{W}(\x,\uu) &\leq -(c_3 - \alpha_1c_4)\norm{\x - h(\uu)}^2
                                 + \alpha_1 c_5\norm{\x-h(\uu)}\norm{\nabla h(\uu)G(h(\uu))}
                                 \notag\\
                               & \hspace{.4cm}
                                 + \alpha_1c_5c_6  \norm{\x-h(\uu)} \norm{\xi},\label{eq:W_minus_lemma}
    \end{align}
  \end{subequations}
  for some $c_1, c_2, c_3, c_4, c_5, c_6 > 0$.
\end{lemma}
The proof is given in Appendix~\ref{sec:proof_lemma_W}.
When $\alpha_1 = 0$, i.e., by assuming a perfect time scale separation between the dynamics of $x$ and $u$, the results in~\eqref{eq:W_conditions_lemma} guarantee $h(u)$ being a global exponentially stable equilibrium of~\eqref{eq:xx_dynamics}. 
\begin{lemma}\label{lemma:S}
  There exists $S: \R^m \to \R$ such that, along the trajectories of~\eqref{eq:uu_dynamics}, it holds 
  \begin{align}
    \dot{S}(\uu)
		&\leq -\alpha_1\norm{\nabla \Fsh(\uu)}^2
		+\alpha_1d_1\norm{\nabla  \Fsh(\uu)}\norm{\x-h(\uu)} 
		+ \alpha_1d_2\norm{\nabla  \Fsh(\uu)}\norm{\xi},\label{eq:S_minus_lemma}
  \end{align}
  for some $d_1, d_2 > 0$.
	Further, $S$ is radially unbounded.
\end{lemma}
 The proof is given in Appendix~\ref{sec:proof_lemma_S}. 
 Condition~\eqref{eq:S_minus_lemma} and LaSalle arguments allow us to claim that, if $\xi = 0$ and $x = h(\uu)$, the trajectory of~\eqref{eq:uu_dynamics} asymptotically enters the set 
$\{u \in \R^m \mid \nabla \Fsh(\uu) = 0\}$.
\begin{lemma}\label{lemma:U}
  There exists a function $U: \R^{2(N-1)d} \to \R$ such that, along the trajectories of~\eqref{eq:xi_dynamics}, it holds
  \begin{subequations}\label{eq:U_conditions_lemma}
    \begin{align}
      b_1 \norm{\xi}^2 &\leq U(\xi) \leq b_2\norm{\xi}^2\label{eq:U_bound_first}
      \\
      \dot{U}(\xi) &\leq -\frac{b_3}{\alpha_2}\norm{\xi}^2 + b_4\norm{\xi}\norm{\x - h(\uu)},\label{eq:U_minus_lemma}
    \end{align}
  \end{subequations}
  for some $b_1, b_2, b_3, b_4 > 0$.
\end{lemma}
The proof is given in Appendix~\ref{sec:proof_lemma_U}.
Lemma~\ref{lemma:U} proves that~\eqref{eq:xi_dynamics} is input-to-state exponentially stable, with input $\x-h(\uu)$. 
Moreover, we also note that the first factor in the right-hand side of~\eqref{eq:U_minus_lemma} can be arbitrarily tuned through $\alpha_2$, i.e., by choosing the timescale separation among the dynamics of $\x$, $\xi$, and $\uu$.

\subsection{Proof of Theorem~\ref{th:convergence}}
\label{sec:proof}

Let $V: \R^{n} \times \R^{m} \times \R^{2(N-1)d}$ be defined as
\begin{align}
  V(\x,\uu,\xi) &:= U(\xi) + W(\x,\uu) + S(\uu),\label{eq:V}
\end{align}
with $W$, $S$, and $U$ given in
Lemma~\ref{lemma:W},~\ref{lemma:S}, and~\ref{lemma:U}.
Let
\begin{align*}
  k_2 &:= \frac{d_1 + c_5}{2}, \quad H_1(\alpha_1) := \begin{bmatrix}
    c_3 - \alpha_1c_4& -\alpha_1k_2\\
    -\alpha_1k_2& \alpha_1\end{bmatrix}.
\end{align*}
Then, by evaluating $\dot{V}(\x,\uu,\xi)$ along the trajectories of~\eqref{eq:continuous_aggregative_gloabal_man} and by using~\eqref{eq:W_minus_lemma},~\eqref{eq:S_minus_lemma}, and~\eqref{eq:U_minus_lemma}, we get
\begin{align}
  \dot{V}(\x,\uu,\xi)
  &\leq
  -\begin{bmatrix}
	\norm{\x - h(\uu)}\\ \norm{\nabla \Fsh(\uu)}
   \end{bmatrix}\T 
   H_1(\alpha_1)
   \begin{bmatrix}
		\norm{\x - h(\uu)}
		\\
		\norm{\nabla \Fsh(\uu)}
	\end{bmatrix}
	+ \alpha_1d_2\norm{\nabla \Fsh(\uu)}\norm{\xi}
	-\frac{b_3}{\alpha_2}\norm{\xi}^2 
	\notag\\
	&\hspace{0.4cm}
	+ (b_4 + \alpha_1 c_5c_6)\norm{\xi}\norm{\x - h(\uu)}
	.\label{eq:dotV}
\end{align}
By Sylvester Criterion, $H_1(\alpha_1) > 0$ if and only if
\begin{align}\label{eq:conditions}
  \begin{cases}
    c_3 > \alpha_1 c_4
    \\
    c_3\alpha_1 > \alpha_1^2(k_2^2 + c_4).
  \end{cases}
\end{align}
Let $\bar{\alpha}_1 := \max\left\{c_3/c_4, c_3/(k_2^2 + c_4)\right\}$. 
This bound quantifies the required ``separation" between the timescales of the dynamics $\x$ and $\uu$.
Indeed, with any $\alpha_1 \in (0, \bar{\alpha}_1)$, both conditions~\eqref{eq:conditions} are satisfied allowing us to claim $H_1(\alpha_1) > 0$. 
Let $h_1(\alpha_1) >0$ be its smallest eigenvalue.
Then, for any $\alpha_1 \in (0,\bar{\alpha}_1)$, we bound~\eqref{eq:dotV} as
\begin{align}
  \dot{V}(\x,\uu,\xi) 
  &\leq
  -h_1(\alpha_1)(\norm{\x-h(\uu)}^2 + \norm{\nabla \Fsh(\uu)}^2)
  + \alpha_1d_2\norm{\nabla \Fsh(\uu)}\norm{\xi}
  -\frac{b_3}{\alpha_2}\norm{\xi}^2 
  \notag\\
  &\hspace{0.4cm}
  + (b_4 + \alpha_1 c_5c_6)\norm{\xi}\norm{\x - h(\uu)}.
  \label{eq:Vodt2}
\end{align}
Let us introduce
\begin{align*}
  e(\x,\uu) &:= \col(\x - h(\uu),\nabla \Fsh(\uu)), \quad
  k_3 := \dfrac{\alpha_1d_2 + b_4 + \alpha_1 c_5c_6}{2}, \quad H_2(\alpha_2)  :=  \begin{bmatrix}
    h_1(\alpha_1)& -k_3\\%
    -k_3& \frac{b_3}{\alpha_2}
  \end{bmatrix}.
\end{align*}
Then, we can bound~\eqref{eq:Vodt2} as 
\begin{align}
  \dot{V}(\x,\uu,\xi) 
  &\leq
    -\begin{bmatrix}
      \norm{e(\x,\uu)}\\ \xi
    \end{bmatrix}\T H_2(\alpha_2)\begin{bmatrix}
      \norm{e(\x,\uu)}\\
      \xi
    \end{bmatrix}.
  \label{eq:dotV5}
\end{align}
This result allows us to quantify the required ``separation'' among the timescales of the dynamics of $\uu$ and $(\w,\z)$.
Let $\bar{\alpha}_2 := \alpha_1b_3h_1(\alpha_1)/k_3^2$.
Then, by Sylvester Criterion, for any $\alpha_2 \in (0,\bar{\alpha}_2)$, it holds $H_2(\alpha_2) > 0$.
Hence, by denoting with $h_2(\alpha_2) > 0$ the smallest eigenvalue of $H_2(\alpha_2)$, the inequality~\eqref{eq:dotV5} leads to 
\begin{align}\label{eq:Vdot_negative}	
  \dot{V}(\x,\uu,\xi)  \leq - h_2(\alpha_2)\|\col(\norm{e(\x,\uu)}, \norm{\xi})\|^2.
\end{align}
Let us study the conditions making the right-hand side of~\eqref{eq:Vdot_negative} is zero. 
To this end, let $\mathcal{U} := \{\uu \in \R^m \mid \nabla \Fsh(\uu) = 0\}$ and
\begin{align}
  E &:= \{(\x,\uu,\xi) \in \R^{n_E} \mid \x = h(\uu),  \uu \in \mathcal{U}, \xi = 0\}\label{eq:E}.
\end{align}
Then $\dot{V}(\x,\uu,\xi) = 0$ for all $(\x,\uu,\xi) \in E$. 
By studying system~\eqref{eq:continuous_aggregative_gloabal_man} restricted to the subset $E$,
we claim that the largest invariant set contained in $E$ for~\eqref{eq:continuous_aggregative_gloabal_man} coincides with $E$ itself. 
Thus, by using the LaSalle Invariance Principle (cf.~\cite[Theorem~4.4]{khalil2002nonlinear}), it holds 
\begin{align}\label{eq:limit}
  \lim_{t\to\infty} \norm{\begin{bmatrix}\x(t)\T&
      \uu(t)\T&
      \xi(t)\T
    \end{bmatrix}\T}_E = 0.
\end{align}
The proof of~\eqref{eq:limit_statement} follows by expressing~\eqref{eq:limit} in the coordinates $\w$ and $\z$ and noting that for all $(\bar{\x},\bar{\uu}) \in \R^n \times \R^m$ such that (i) $\bar{\x} = h(\bar{\uu})$ and (ii) $\nabla \Fsh(\bar{\uu}) = 0$, it holds $(\bar{\x},\bar{\uu}) \in X$~\cite[Prop.~3.1]{hauswirth2020timescale}.
As for the second claim, we pick any $\bar{\uu} \in \R^{m}$ being an isolated stationary point and a local minimum of $\Fsh$.
Thus, by definition, there exists a neighborhood $\bar{\mathcal{U}} \subseteq \R^{m} \setminus \{\bar{\uu}\}$ of $\bar{\uu}$ such that $\Fsh(\uu) > \Fsh(\bar{\uu})$ and $\nabla \Fsh (\uu) \ne 0$ for all $\uu \in \bar{\mathcal{U}}$.
By using these facts and looking at the definition of $V$ (cf.~\eqref{eq:V}) and the inequality about $\dot{V}$ given in~\eqref{eq:Vdot_negative}, we guarantee that there exists a neighborhood $\bar{\cS} \subseteq \R^{n_E} \setminus \{\col(h(\bar{\uu}),\bar{\uu},0)\}$ of $\col(h(\bar{\uu}),\bar{\uu},0)$ such that $V(\x,\uu,\xi) > V(h(\bar{\uu}),\bar{\uu},0)$ and $\dot{V}(\x,\uu,\xi) < 0$ for all $(\x,\uu,\xi) \in \bar{\cS}$.
Therefore, the proof follows by using, e.g.,~\cite[Th.~4.1]{khalil2002nonlinear}.

\section{Multi-robot Surveillance}
\label{sec:numerical_simulations}

\begin{figure}
	\centering
	\subfloat[][Optimality error.]
	{\includegraphics[scale=.9]{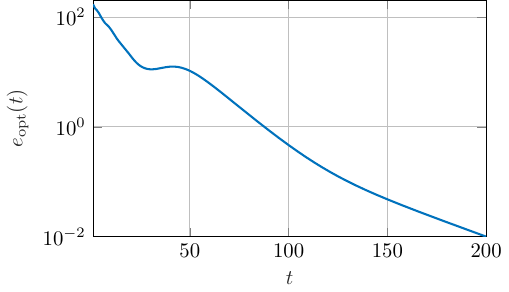}}\quad
	\subfloat[][Error of $\w$ and $\z$.]
	{\includegraphics[scale=.9]{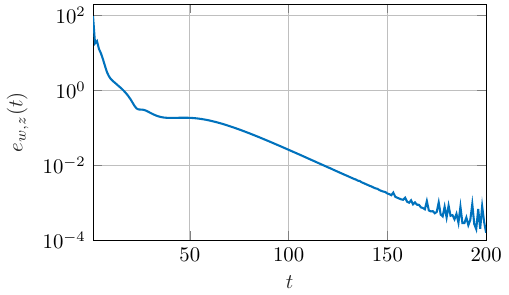}}
	\caption{Multi-robot surveillance: errors evolution.}
	\label{fig:error}
\end{figure}

\begin{figure}
	\centering
	\subfloat[][Initial configuration.]
	{\includegraphics[scale=.4]{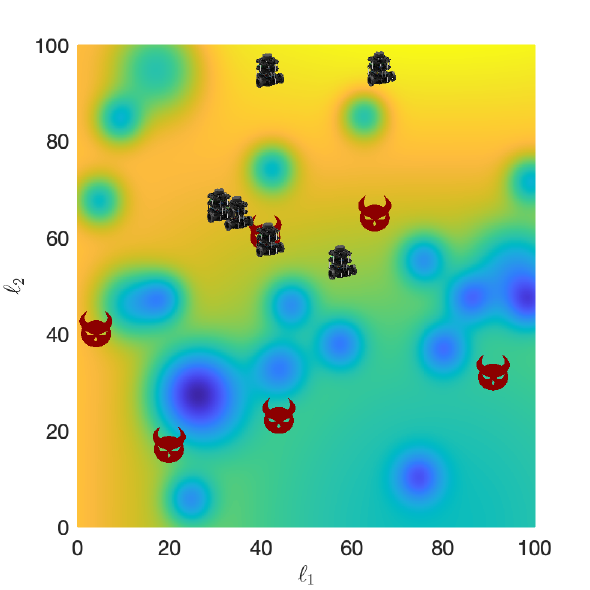}}\quad
	\subfloat[][Final configuration.]
	{\includegraphics[scale=.4]{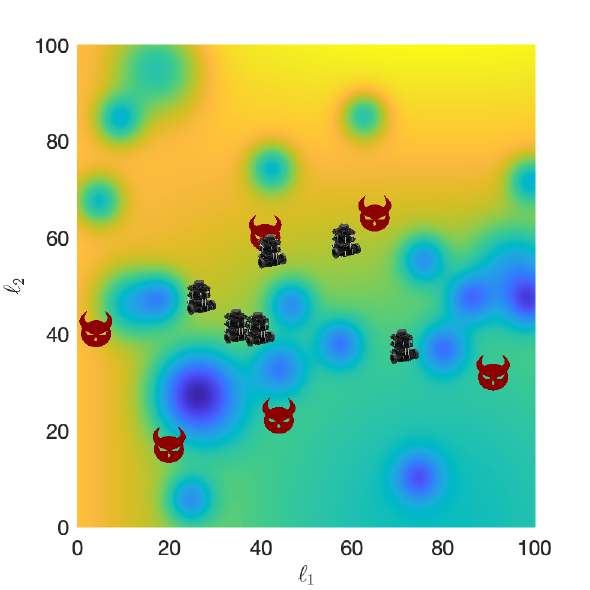}}
	\caption{Multi-robot surveillance: nonconvex scenario.}
	\label{fig:screen}
\end{figure}

\begin{figure}
	\centering
	\subfloat[][Initial configuration.]
	{\includegraphics[scale=.4]{figs/initial_configuration.pdf}}\quad
	\subfloat[][Final configuration.]
	{\includegraphics[scale=.4]{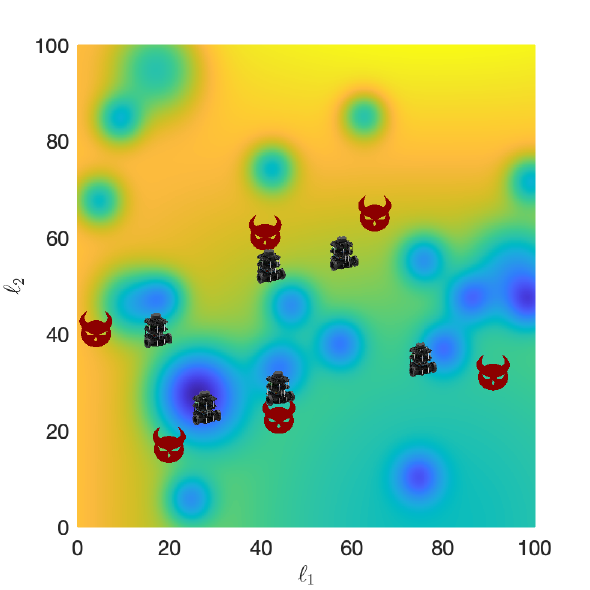}}
	\caption{Multi-robot surveillance: strongly convex scenario.}
	\label{fig:screen_sc}
\end{figure}

\begin{figure}
	\centering
	\subfloat[][Optimality error.]
	{\includegraphics[scale=.9]{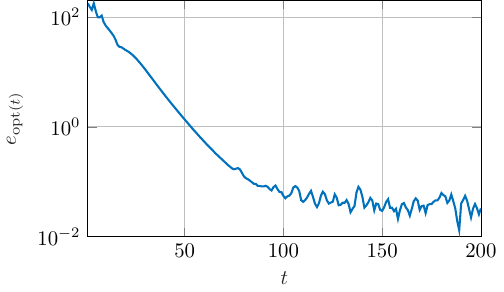}}\quad
	\subfloat[][Error of $\w$ and $\z$.]
	{\includegraphics[scale=.9]{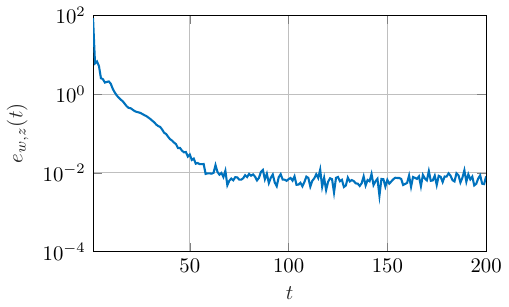}}
	\caption{Multi-robot surveillance: disturbed case.}
	\label{fig:error_dist}
\end{figure}

In this section, we use \algo/ to address a multi-robot surveillance scenario. %
We consider a network of $N$ mobile robots, whose planar position is $\x_i \in \R^2$, that aim to
surveil a collection of $N$ intruders each
located at $s_i \in \R^2$. 
Given the orientation $\theta_i \in \R$ of robot $i$, we describe its dynamics as %
\begin{align}\label{eq:unicycle}
 \dot{\x}_i &= \begin{bmatrix}
 \cos(\theta_i)
 \\
 \sin(\theta_i)
 \end{bmatrix}v_i
 ,\quad \dot{\theta}_i = \omega_i,
\end{align} 
where $v_i, \omega_i  \in  \R$ are low-level inputs denoting linear and angular speed, respectively. 
Let $\uui \in \R^2$ be a reference position, then~\cite{terpin2022distributed} proposes the low-level control
\begin{subequations}\label{eq:low_level}
\begin{align}
 v_i(\x_i,\theta_i,u_i) &= k_i\norm{\x_i - u_i}\cos(\tilde{\theta}_i(\x_i,\theta_i))
 \\
 \omega_i(\x_i,\theta_i,u_i) &= \tfrac{k_i}{\norm{\x_i - u_i}}\cos(\tilde{\theta}_i(\x_i,\theta_i))\sin(\tilde{\theta}_i(\x_i,\theta_i))
 +\tfrac{k_i}{\norm{\x_i - u_i}}\sin(\tilde{\theta}_i(\x_i,\theta_i)),
\end{align}
\end{subequations}
with $k_i  >  0$, $\tilde{\theta}_i(\x_i,\theta_i)  =  \text{atan}2(\x_{i,1},\x_{i,2})  -  \theta_i$, and $\x_i  :=  \col(\x_{i,1},\x_{i,2})$. 
Thus, the closed-loop dynamics reads as
\begin{subequations}\label{eq:unicycle_closed_loop}
\begin{align}
 \dot{\x}_i &= \begin{bmatrix}
 \cos(\theta_i)
 \\
 \sin(\theta_i)
 \end{bmatrix}v_i(\x_i,\theta_i,\uu_i)
 \\
 \dot{\theta}_i &= \omega_i(\x_i,\theta_i,\uu_i).
\end{align} 
\end{subequations}
The point $\col(\uu_i,0)$ is an almost globally exponentially
stable equilibrium for~\eqref{eq:unicycle_closed_loop} for all $\uu_i$~\cite[Lemma~2.1]{terpin2022distributed}. 
\begin{remark}
	We note that Assumption~\ref{ass:steady_state} requires global exponential stability for the entire plant state, while~\eqref{eq:unicycle_closed_loop} guarantees this only for the portion $\x_i$ of the state. %
	However, as shown in~\cite{terpin2022distributed}, one may easily modify this assumption to handle the more general case in which it can be only guaranteed exponential stability for a portion of the state. 
	In this case, the cost must only depend on the exponentially stable portion of the state. \oprocend
\end{remark}
As for the environment, we consider a nonconvex scenario in which altitude changes and $n_c \in \mathbb{N}$ crevasses are present.
Let $\col(\ell_1,\ell_2)$ be the planar coordinates describing a given location. 
Then, we model the altitude profile of the environment through a function
$\alt: \R^2 \to \R$ given by the sum of a sinusoidal term and a series of Gaussian functions modeling the crevasses, namely
\begin{align}\label{eq:altitude}
 &\alt(\ell_1,\ell_2) = - a_1\cos(\rho \ell_1)\sin(\rho \ell_2) 
 - \sum_{g=1}^{n_c} a_{c,g} \exp\big(-\tfrac{1}{s_g}(\ell_1 - \mu_{g,1})^2 + (\ell_2 - \mu_{g,2})^2\big),
\end{align}
where $a_1, \rho > 0$ are the amplitude and the frequency of the
sinusoidal term, while the parameters $a_{c,1}, \dots, a_{c,n_c}, s_1, \dots, s_{n_c} > 0$ characterize the Gaussian functions with centers located in $(\mu_{1,1},\mu_{1,2}), \dots, (\mu_{n_c,1},\mu_{n_c,2})$.
This environment profile gives rise to a nonconvex problem.
The surveillance strategy of the team consists of a trade-off between the
following competing objectives: each robot (i) tries to stay close to the
intruder, (ii) tries to occupy locations with higher altitudes, and (iii) tries
to stay close to the weighted center of mass.
This scenario falls into the aggregative framework by setting each cost function $f_i$ as
\begin{align}
 f_i(\x_i,\sigma(\x)) &= \gamma_1\norm{\x_i - s_{i}}^2 - \alt(\x_{i,1},\x_{i,2})
 + \gamma_2\norm{\x_i - \sigma(\x)}^2,\label{eq:surveillance}
 \end{align} 
 where $\gamma_1, \gamma_2 > 0$, while the term $-\alt(\x_i)$ increases the cost according to the altitude of the location $\x_i$ (cf.~\eqref{eq:altitude}).
 Further, we choose $\sigma(x)$ as the weighted center of mass of the defending team, namely 
 $\sigma(\x) = \frac{1}{N}\sum_{i=1}^{N}\beta_{i}\x_i$,%
 for some weights $\beta_{i} > 0$. 
 In detail, we consider a network of $N= 6$ agents communicating according to an \er/ graph with connectivity parameter $p = 0.4$.
Further, we set $\gamma_1 = 1$, $\gamma_2 = 0.3$, $n_c = 5$, and randomly generate the weights $\beta_i \in (0,1)$, the amplitudes $a_{c,1}, \dots, a_{c,n_c} \in [0,5]$, the terms $s_1, \dots, s_{n_c} \in (5,10)$, and the locations $\mu_1 := \col(\mu_{1,1},\mu_{1,2}), \dots, \mu_{n_c} := \col(\mu_{n_c,1},\mu_{n_c,2})$, $y_1, \dots, y_N$, and $b$ within the interval $[0,100]^2$. 
As for the sinusoidal terms, we choose $a_1 = 10$ and $\rho = 0.02$.
As for the algorithm parameters, we empirically tuned $\alpha_1 = \alpha_2 = 0.75$ (with $\alpha_1 = 7$ and $\alpha_2 = 1$ the convergence properties are lost), while $\x_{i}(0)$ and $\uu_i(0)$ are randomly selected. 
The simulations are carried out by using the function ode45 of Matlab.
Let $e_{\text{opt}} := \norm{\col(\x(t) - \uu(t),\nabla \Fsh(\uu(t)))}$ and $e_{\w\z}(t) := \norm{\col(\w(t) - \pi^{\w}(x(t)),\z(t) - \pi^{\z}(x(t)))}$.
 As predicted by Theorem~\ref{th:convergence}, Fig.~\ref{fig:error} shows that $e_{\text{opt}}$ and $e_{\w\z}$ asymptotically converges to $0$. %
Considering the same simulation, Fig.~\ref{fig:screen} provides the initial and final configuration of the team.
Each robot icon denotes an agent of the surveillance team, while each devil icon denotes an intruder. 
The color of the background represents the altitude: blue background denotes the lowest locations, while yellow background denotes the highest ones.
Fig.~\ref{fig:screen} highlights the role played by the altitude in determining the final configuration achieved by the agents. 
Indeed, some of the robots remain far from their intruders because closer locations would have lower altitudes.
To emphasize this point, we repeat the simulations neglecting the altitude $\alt$ in the cost, i.e., by setting $\gamma_2 = 0$ into~\eqref{eq:surveillance}.
Fig.~\ref{fig:screen_sc} provides the initial and final team configuration of such a simulation. %
Here, differently from the previous case, the robots go closer to their associated intruders thus occupying locations with low altitudes.
In both cases, the robots arrange themselves inside the polygon whose vertices coincide with the positions occupied by the intruder.
In fact, the outer configurations suffer a higher cost due to the term $\norm{\x_i - \sigma(\x)}^2$.
Finally, we repeat the simulations with disturbances $d_i \in \R^{n_i}$ acting on each plant dynamics~\eqref{eq:local_plant}, namely $\dot{\x}_i = p_i(\x_i,\uu_i) + d_i$.
Fig.~\ref{fig:error_dist} reports the evolution of $e_{\text{opt}}$ and $e_{\w,\z}$ by using the same data of the nominal simulation and $d_i$ randomly extracted from the interval $[-0.5,0.5]$ with uniform distribution for each time $t$.
In the presence of disturbances, the errors $e_{\text{opt}}$ and $e_{\w,\z}$ do not exactly converge to zero but remain bounded thus showing the typical behavior of perturbed systems in the neighborhood of asymptotically stable points.

\section{Conclusion}

We proposed \algo/, i.e., a novel distributed feedback optimization law for the aggregative framework and nonlinear dynamics. 
Our scheme relies on the communication among the agents and the measurements about the optimization problem information. 
With tools from system theory, we proved that our method steers the network to a stationary point of the optimization problem. 
Some simulations on a surveillance scenario confirmed these results. %

\vspace{-.1cm}

\appendix
\numberwithin{equation}{section}

\renewcommand{\thesection}{\Alph{section}}

\section{Proof of Lemma~\ref{lemma:W}}
\label{sec:proof_lemma_W}
By using the Converse Lyapunov Theorem (cf.~\cite[Theorem~5.17]{sastry2013nonlinear}), the exponential stability of $h(\uu)$, and the Lipschitz continuity of $h$ (cf. Assumption~\ref{ass:steady_state}), there exists $W: \R^n \times \R^m \to \R$ and $c_1, c_2, c_3, c_5 > 0$ such that
			\begin{subequations}\label{eq:W_conditions}
				\begin{align}
					c_1 \norm{\x - h(\uu)}^2 \leq W(\x,\uu) &\leq c_2\norm{\x - h(\uu)}^2\label{eq:W_first}
					\\
					\nabla_1 W(\x,\uu)\,p(\x,\uu) &\leq -c_3\norm{\x - h(\uu)}^2\label{eq:W_minus}
					\\
					\nabla_2 W(\x,\uu) &\leq c_5\norm{\x - h(\uu)},\label{eq:W_third}
				\end{align}
			\end{subequations}
			Since~\eqref{eq:W_first} matches~\eqref{eq:W_upepr_bound}, we need only to show~\eqref{eq:W_minus_lemma}. 
			Along the trajectories of~\eqref{eq:xx_dynamics} and~\eqref{eq:uu_dynamics}, it holds
			\begin{align}
				&\dot{W}(\x,\uu) =
				\nabla_1 W(\x,\uu)p(\x,\uu) +\nabla_2 W(\x,\uu)\dot{\uu}
				\notag\\
				&
				\stackrel{(a)}{\leq}
				-c_3\norm{\x-h(\uu)}^2 +\nabla_2 W(\x,\uu)\dot{\uu}
				\notag\\
				&\stackrel{(b)}{\leq}
				-c_3\norm{\x-h(\uu)}^2 + c_5\norm{\x-h(\uu)}\norm{\dot{\uu}},\label{eq:dotW_intemrediate}
			\end{align}
			where in $(a)$ we use~\eqref{eq:W_minus}, and in $(b)$ we use the Cauchy-Schwarz inequality with condition~\eqref{eq:W_third}.
			Since $G(\x) = 
				\Go(\x,\1\sigma(\x)) 
				+\nabla \phi(\x)\1\1\T\Gt(\x,\1\sigma(\x))/N$, if we add and subtract $\alpha_1\nabla h(\uu)\Go(\x,\1\sigma(\x))$ into~\eqref{eq:uu_dynamics}, we get
			\begin{align}
				\dotuu &= -\alpha_1\nabla h(\uu)G(\x) - \alpha_1\nabla\phi(\x)\OR\xi
				-\alpha_1\nabla h(\uu)\big(\Go(\x,\RO\xi + \1\sigma(\x)) -\Go(\x,\1\sigma(\x))
				\notag\\
				&\hspace{0.4cm}
				+\nabla \phi(\x)\tGt(\x,\xi)\big).\label{eq:uu_dynamics_changed}
			\end{align}
			Moreover, by using the Lipschitz continuity properties given in Assumption~\ref{ass:steady_state} and~\ref{ass:lipschitz}, we can write
			\begin{subequations}\label{eq:inqualities}
			\begin{align}
				\|\Go(\x,\RO\xi + \1\sigma(\x)) -\Go(\x,\1\sigma(\x))\| &\leq  L_1 \|\xi\| 
				\\
				\|\tGt(\x,\xi)\| &\leq L_2\|\xi\|
				\\
				\|\nabla \phi(\x)\| &\leq  L_3
				\\
				\|\nabla h(\uu)\| &\le  L_h.		
			\end{align}
			\end{subequations}
			Then, we combine~\eqref{eq:uu_dynamics_changed}, the Cauchy-Schwarz inequality, and the bounds~\eqref{eq:inqualities} to obtain
			\begin{align}
					\norm{\dotuu} &\leq \alpha_1\norm{\nabla h(\uu)G(\x)} 
					+ \alpha_1L_h (L_1 + (1+L_2)L_3) \|\xi\|
					\notag\\
					&\stackrel{(a)}{\leq}
					\alpha_1\norm{\nabla h(\uu)G(h(\uu))} 
					+ \alpha_1\norm{\nabla h(\uu)G(\x)  - \nabla h(\uu)G(h(\uu))}
					+ \alpha_1L_h (L_1 + (1+L_2)L_3)\|\xi\|
					\notag\\
					&\stackrel{(b)}{\leq}
					\alpha_1 \norm{\nabla h(\uu)G(h(\uu))} 
					+ \alpha_1L_hL_0\norm{x  -  h(\uu)}
					+ \alpha_1L_h(L_1 + (1+L_2)L_3)\|\xi\|,\label{eq:uu_dynamics_perturbation}
			\end{align}
			where in $(a)$ we add $\pm\nabla h(\uu)G(h(\uu))$ and use the triangle inequality, while $(b)$ uses the Lipschitz continuity of $h$ and $\nabla F$ (cf. Assumptions~\ref{ass:steady_state} and~\ref{ass:lipschitz}).
			The proof follows using $\nabla h(\uu)G(h(\uu)) = \nabla \Fsh(\uu)$ and~\eqref{eq:uu_dynamics_perturbation} into~\eqref{eq:dotW_intemrediate}. %

\section{Proof of Lemma~\ref{lemma:S}}
\label{sec:proof_lemma_S}

We set $S(\uu) := \Fsh(\uu)$.
			Then, $S$ is radially unbounded (cf. Assumption~\ref{ass:lipschitz}).
			Along the orbits of~\eqref{eq:uu_dynamics}, it holds
		   \begin{align}
			\dot{S}(\uu) 
			&= (\nabla \Fsh(\uu))\T\dot{\uu}
			 \notag\\
			&
			\stackrel{(a)}{\leq}
			 -\alpha_1(\nabla \Fsh(\uu))\T(\nabla h(\uu)G(\x)) 
			+ \alpha_1d_2\norm{\nabla\Fsh(\uu)}\norm{\xi}
			\notag\\
			&\stackrel{(b)}{=}
			 -\alpha_1\norm{\nabla \Fsh(\uu)}^2
			 + \alpha_1d_2\norm{\nabla \Fsh(\uu)}\norm{\xi}
			 -\alpha_1(\nabla \Fsh(\uu))\T\big(\nabla h(\uu)G(\x)
			 - \nabla h(\uu)G(h(\uu))\big) 
			\notag\\
			&
			 \stackrel{(c)}{=}
			 -\alpha_1\norm{\nabla\Fsh(\uu)}^2
			 + \alpha_1d_2\norm{\nabla \Fsh(\uu)}\norm{\xi}
			 +\alpha_1L_hL_0\norm{\nabla \Fsh(\uu)}\norm{\x-h(\uu)},\notag%
		\end{align}
		where $(a)$ uses~\eqref{eq:uu_dynamics_changed} and~\eqref{eq:inqualities} setting $d_2 = L_h(L_1 + (1+L_2)L_3)$, in $(b)$ we add $\pm\nabla\Fsh(\uu) = \nabla h(\uu)G(h(\uu))$, and $(c)$ uses the Lipschitz continuity of $h$ and $G$ (cf. Assumptions~\ref{ass:steady_state},~\ref{ass:lipschitz}) and the Cauchy-Schwarz inequality.
		The proof follows by setting $d_1 = L_hL_0$.

\section{Proof of Lemma~\ref{lemma:U}}
\label{sec:proof_lemma_U}

In light of Assumption~\ref{ass:network}, the matrix $-R\T LR$ is Hurwitz. 
Thus, given $q_1, q_2 > 0$, there exist $P_1, P_2 \in \R^{(N - 1)d \times (N - 1)d}$ such that
		\begin{subequations}
			\label{eq:lyapunov_equations}
		\begin{align}
			-P_1R\T LR - (R\T LR)\T P_1 &= -q_1 I
			\\
			-P_2R\T LR - (R\T LR)\T P_2 &= -q_2 I,\label{eq:LjapQ2}
		\end{align}
		\end{subequations}
		and $P_1 = P_1\T > 0$, $P_2 = P_2\T > 0$.
		Then, let
		$U(\xi) := \xi\T P \xi$,
		with $P := \blkdiag(P_1,P_2)$.
		Hence, the conditions~\eqref{eq:U_bound_first} are verified by using the eigenvalues of $P$. %
		In order to show~\eqref{eq:U_minus_lemma}, let $\xi_1,\xi_2 \in \mathbb{R}^{(N-1)d}$ be such that $\xi = \col(\xi_1,\xi_2)$.
		Then, by using~\eqref{eq:xi_dynamics} and~\eqref{eq:lyapunov_equations}, we get
				\begin{align}
					\dot{U}(\xi) &=
					 - \frac{q_1}{\alpha_2}\norm{\xi_1}^2
					 - \frac{q_2}{\alpha_2}\norm{\xi_2}^2
					+ \frac{2}{\alpha_2}\xi_2\T
						P_2R\T L\tGt(\x,\xi)
					- 2\xi\T P\,\nabla\bar{\psi}(\x)\, p(\x,\uu)
					\notag
					\\
					&\stackrel{(a)}{\leq}
					 - \frac{q_1}{\alpha_2}\norm{\xi_1}^2- \frac{q_2}{\alpha_2}\norm{\xi_2}^2
					+ \frac{2L_2}{\alpha_2}\norm{P_2R\T L}\norm{\xi_2}\norm{\xi_1}
					- 2\xi\T P\,\nabla\bar{\psi}(\x)\, p(\x,\uu),\label{eq:dotU}
			\end{align}
				where in $(a)$ we use the Cauchy-Schwarz inequality and the Lipschitz continuity of $\nabla_2 f_i$ (cf. Assumption~\ref{ass:lipschitz}).
				Now, we fix $q_2 > 0$ and compute $P_2$ such that~\eqref{eq:LjapQ2}.
				Let
				\begin{align*}
					k_1(q_2) := L_2\norm{P_2R\T L}, \quad 
				\tilde{Q} := \begin{bmatrix}
					q_1& -k_1(q_2)\\
					-k_1(q_2)& q_2
				\end{bmatrix}.
			\end{align*}
				Then, we rewrite~\eqref{eq:dotU} as
				\begin{align}
					\dot{U}(\xi)
						 &\leq 
						  - \frac{1}{\alpha_2}\begin{bmatrix}
							\norm{\xi_1}\\
							\norm{\xi_2}
						\end{bmatrix}\T \tilde{Q}  \begin{bmatrix}
							\norm{\xi_1}\\
							\norm{\xi_2}
						\end{bmatrix}
							  -  2\xi\T P \nabla\bar{\psi}(\x)\, p(\x,\uu).\notag
				\end{align}
				Let us choose $b_3 \in (0, q_2)$ and $q_1 > \tfrac{b_3(q_2 - b_3)+k_1(q_2)^2}{q_2-b_3}$. 
				Then, it holds $\tilde{Q} > b_3 I$ which allows us to write
				\begin{align}
					\dot{U}(\xi)&\leq
					-\frac{b_3}{\alpha_2}\norm{\xi}^2  + 2\xi\T P \nabla\bar{\psi}(\x)\, p(\x,\uu).\label{eq:dotU2}
				\end{align}
				Since $p(h(\uu),\uu) = 0$ (see Assumption~\ref{ass:steady_state}), it holds
				\begin{align}
					\xi\T  P \nabla\bar{\psi}(\x) p(\x,\uu) &=  \xi\T P\nabla\bar{\psi}(\x)(p(\x,\uu) - p(h(\uu),\uu))%
					\notag\\
					&\stackrel{(a)}{\leq} %
					L_p \norm{P}\norm{\nabla \bar{\psi}(\x)}\norm{\xi}\|x - h(\uu)\|,\label{eq:bound_proof2}
				\end{align}
				where in $(a)$ we use Assumption~\ref{ass:steady_state}.
				By using the Cauchy-Schwarz inequality, Assumption~\ref{ass:lipschitz}, $\|R\| = 1$, and $\|1_N 1_n^\top \| = \sqrt{Nn}$, we get $\|\nabla\bar{\psi}(\x)\|  
		\le (L_2\sqrt{Nn}+L_3)$.%
		The proof follows by  combining the latter with~\eqref{eq:dotU2} and~\eqref{eq:bound_proof2} and setting $b_4 := \tfrac{L_p\|P\|(L_2\sqrt{Nn} + L_3)}{2}$.

\end{document}